\documentclass{amsart}
\usepackage[margin=1in]{geometry}
\usepackage{amssymb, amscd, amsmath, url}

\newtheorem{theorem}{Theorem}[section]

\newtheorem{lemma}{Lemma}[section]

\newtheorem{remark}{Remark}[section]

\numberwithin{equation}{section}

\usepackage{amsmath, amsthm}

\usepackage[linesnumbered,ruled,vlined]{algorithm2e}
\SetCommentSty{mycommfont}
\SetKwInput{KwInput}{Input}                
\SetKwInput{KwOutput}{Output}
\usepackage{graphicx} 
 
\usepackage[version=4]{mhchem}
\usepackage{stmaryrd}
\usepackage{multirow}

\usepackage{xcolor}
\usepackage[linesnumbered,ruled,vlined]{algorithm2e}

\SetCommentSty{mycommfont}

\begin{document}

\title[]{A probabilistic representation for the gradient in a linear parabolic PDE with Neumann boundary condition}



\author{\bfseries Abdelatif Bench\'{e}rif Madani $^1$}
\address{$^1$ Laboratory of Fundamental and Numerical Mathematics (LMFN), Department of mathematics, University Setif - 1 - Ferhat Abbas, Algeria.}
\maketitle

\begin{abstract}We give a probabilistic representation for the gradient of a 2nd order linear parabolic PDE $
\partial_{t}u(t,x)=(1/2)a^{ij}\partial_{ij}u(t,x)+b^{i}\partial_{i}u(t,x)$ with Cauchy initial condition $u(0,x)=f(x)$ and Neumann boundary condition in a (closed) convex bounded smooth domain $D$ in $\mathbb{R}^{d}$, $d\geq 1$. The idea is to start from a penalized version of the associated reflecting diffusion $X^{x}$, proceed with a pathwise derivative, show that the resulting family of $\nu$-directional Jacobians is tight in the Jakubowski S-topology with limit $J^{x,\nu}$, solution of a certain linear SDE, and set $\mathbb{E}\left(\nabla f(X^{x}(t))\cdot J^{x,e_{i}}(t)\right)$ for the gradient $\partial_{i}u(t,x)$, where $x\in D$, $t\geq 0$, $e_{i}$ the canonical basis of $\mathbb{R}^{d}$ and $f$, the initial condition of the semigroup of $X^{x}$, is differentiable. Some more extensions and applications are discussed in the concluding remarks.

\end{abstract}

\section{Introduction and statement of the result}

In this paper we study the pointwise differentiation of the flow of the reflecting It\^{o} diffusion $X^{x}$ in a closed domain $D\subset \mathbb{R}^{d}$ in the direction of the unit inward pointing normal field $\gamma$ at the boundary $\partial D$
\begin{equation}\label{Xdiffus} 
X^{x}_{t}=x+\int_{0}^{t}b(X^{x}_{s})ds+\int_{0}^{t}\sigma(X^{x}_{s})dB_{s}
+\int_{0}^{t}\gamma(X^{x}_{s})dL^{x}_{s},
\end{equation} 
where $L(t)$ is the boundary local time and the coefficient $b$ (resp. $\sigma$) is a vector function (resp. a $d\times d$-matrix function whose $(i,j)$-th entry is $\sigma_{j}^{i}$) in $\mathcal{C}_{b}^{1}(\mathbb{R}^{d})$. Moreover, we take $a=\sigma\sigma^{*}$ to be uniformly non degenerate: for some $c>0$, we have $(a(x)(\eta),\eta)\geq c$ for all $x\in D$ and $\eta\in\mathbb{R}^{d}$. We assume we are working in a probability space $(\Omega,\mathcal{F},\mathcal{F}_{t},\mathbb{P})$ satisfying the usual conditions. The process $L$ is constant over inward excursion times and grows only on the (time) null set $Z^{x}=\lbrace t|X^{x}_{t}\in \partial D\rbrace$ a.s.. The main motivation for our work is the study of the gradient of the following Cauchy second order linear parabolic boundary value problem 
\begin{equation}\label{EqParbol}
\begin{aligned}
& \partial _{t}u=\dfrac{1}{2}\sum_{i,j=1}^{d}a^{ij}(x)\partial _{x^{i}x^{j}}^{2}u + \sum_{i=1}^{d}b^{i}(x)\partial_{x^{i}}u, \\
& u(0,x)=f(x),\: \partial _{\gamma}u =0 \: on \:\partial D, 
\end{aligned} 
\end{equation}
where $\partial _{\gamma}$ stands for the directional derivative along the field $\gamma$ at $\partial D$. The solution $u(t,x)$ of \eqref{EqParbol} has the probabilistic representation $u(t,x)=\mathbb{E}f(X^{x}(t))$. Taking formal derivatives, that is making attempts at studying the gradient system $\partial _{i}u(t,x)$ ($i=1,...,d$), on the flow $x\mapsto X^{x}(t)$ associated with the diffusion $X$ in the above probabilistic formula for the solution $u$, we should have for any direction $\nu$ as $\epsilon\rightarrow 0$
\begin{equation}\label{Ratio} 
\partial_{\nu}u(t,x)=\mathbb{E}\left(\partial f(X^{x}(t))\cdot\lim \dfrac{X^{x+\epsilon\nu}(t)-X^{x}(t)}{\epsilon}\right).
\end{equation}
In the case without boundary, it is well known \cite{Kun} that the flow is invertible and can be made arbitrarily regular depending on the regularity of the coefficients $a,b$ and the above manipulation is indeed rigorous for differentiable $f$. The Jacobian matrix $J^{x}(t)=[\partial_{ij}^{2}X^{x}(t)]$ satisfies the linear equation (below $\mathbb{I}$ is the $d\times d$ identity matrix)
\begin{equation}\label{JacobianFree} 
J^{x,\nu}_{t}=\mathbb{I}+\int_{0}^{t}\partial b(X_{s}^{x})J^{x,\nu}_{s}ds+\int_{0}^{t}\partial \sigma_{ j}(X^{x}_{s})J^{x,\nu}_{s}dB_{s}^{j}
\end{equation}
where $\partial b(x)$ (respectively $\partial \sigma_{ \cdot}(x)$) is the matrix $\partial_{x^{j}} b^{i}(x)$ (respectively the matrix $\partial_{x^{j}} \sigma_{ \cdot}^{i}(x)$).  The general case \eqref{EqParbol} with boundary is drastically different since the flow is no more invertible and indeed discontinuous. Moreover, for systems of PDEs, the boundary data is quite varied and has to be addressed; as a rule, it seems to be more naturally expressed within a geometric framework, see e.g. \cite{Spi} and \cite{MitTay}. In order to gain some insights and put our contribution into perspective, it is worth describing in more detail the existing works about the subject. 

The whole story seems to have begun with \cite{Air}. She studied a probabilistic representation for a system of harmonic functions $(f_{1},...,f_{d})$ with the boundary condition $(N(x)-N^{\bot}(x)\partial_{\gamma})f(x)=g(x)$  with possibly degenerate matrices $N,N^{\bot}$, $g$ is a function. Her motivation, which is of interest in many areas of mathematics, was clearly geometric since this kind of boundary behaviour is reminiscent of the scalarization (on the orthonormal frame bundle) of the absolute boundary condition for the heat equation on 1-forms in Riemann manifolds, see Section 4 of \cite{Hsu2}. The central phenomenon, arrived at using the method of singular perturbations, was to pass to the limit in the following process: (keeping essentially Airault's original notation) $2_{\epsilon}(t)=exp \star(-\int_{0}^{t}(N(X_{s})/\epsilon)dL(s))$ where $exp \star (\cdot)$ refers to the multiplicative integral. The limit exists, loosely speaking, because of an interesting averaging principle (a well known paradigm in the homogenization and ergodic theory communities e.g.) : the blow up as $\epsilon\rightarrow 0$ is compensated with the local time being (time) a.a. constant. The scalar case $d=1$ is well understood, for an updated analytical account see \cite{Tai}. The paper \cite{Air} was taken up in the hard going chapter 5, Section 6, of \cite{IkeWat} using again excursion theory and it was much amplified in \cite{Hsu2}: the problem was clearly identified as giving a probabilistic representation for the solution of the heat equation (for the Hodge-De Rham Laplacian) for differential  1-forms with the absolute boundary condition in a Riemann manifold with boundary. Here, a curved manifold gives rise to a couple of notions of stochastic differentials. For the anecdote, \cite{Hsu2} observes on p. 351 that the manipulations in \cite{IkeWat} "seem to have been created especially for this problem" and departed from the pathwise derivative stand-point, and its associated excursion theory, in favour of a Feynman-Kac argument. In the long and hard going paper \cite{ArnLi}, still in the manifold case, another explanation for \cite{IkeWat} was given in an Appendix therein. In \cite{ArnLi} we do find a penalization construction of reflected Brownian motion and a subsequent pathwise covariant derivation. However, much like \cite{Pil}, the penalization is so strong that the boundary local time is altogether avoided and the reflected diffusion is wholly obtained from the interior of $D$.

In $\mathbb{R}^{+}=[0,\infty)$, directly differentiating the flow is carried out in \cite{BosCisTal}, (see also \cite{ManMas}): setting $J^{x}(t)=dX^{x}(t)/dx$ and letting $\zeta$ be the hitting time of the origin, it turns out that $dJ^{x}(t)=db/dx(X^{x}(t))J^{x}(t)dt+d\sigma/dx(X^{x}(t))J^{x}(t)dB(t)$ when $t<\zeta$ and $J^{x}$ vanishes identically afterwards, $J^{x}(0)=1$. The substitute for the latter property of the Jacobian in $\mathbb{R}^{d}$ is the projection onto the tangent space phenomenon. The reflecting flow was properly differentiated (i.e. performing pathwise the limit in \eqref{Ratio}) in \cite{DoeZam} (normal reflection in the first orthant and $\sigma$ constant), \cite{And1} (oblique reflection in a convex polyhedron and $\sigma$ constant), \cite{Bur} (normal reflection in arbitrary smooth $D$ but with $\sigma$ constant and $b=0$) (it has been revisited in \cite{BurLee} which seems to be mostly a display of techniques), \cite{And2} (normal reflection in arbitrary smooth $D$ and $\sigma$ constant), \cite{LipRam} (oblique reflection in a convex polyhedron but with general coefficients depending on a parameter and an extra sensitivity analysis).

There are also, when there is too much randomness, indirect differentiations. In \cite{Pil} a general diffusion reflecting in the normal direction in the upper half-space $\mathbb{R}^{d}_{+}$ is considered and an appealing heuristics is provided : at the local minimas of the surface $\mathbb{R}^{d}\ni x\rightarrow X^{x,d}\in \mathbb{R}_{+}$, the last component of $J_{t}^{x,\nu}$ must be re-set to zero, and the actual differentiation is recovered "from the back door", so to speak. Note that \cite{LipRam} is also, in some sense, indirect since a "derivative process" is first (ingenuously) suggested and then proved to be the right one thanks to properties of the derivatives of the so called extended Skrorohod map beginning with an orthant (much of the difficulties concern the multiplicity of the reflection directions) and first applications are in Queuing theory. Observe that no curvature terms occur in \cite{LipRam}, see also Section \ref{Weingrtn} below. Our representation here seems to be in the same vein but first applications are in PDE theory, inspired by \cite{DoeZam} and \cite{And1}.

A crucial ingredient in the foregoing treatments is the fact that the Skorohod map, the extended Skorohod map and subsequently $x\mapsto L^{x}(t)$ are Lipschitz; in the general situation this seems to be, to the best of our knowledge, an open problem. Note that the flow associated with reflecting diffusions has been studied in many more ways, concerning coalescence for example. 

As appears from the account above, all the authors have had to assume some additional assumption or another. Therefore, the pointwise tracking in $x,t,\epsilon$ the ratio in \eqref{Ratio} in the general case (i.e. general coefficients and a non-convex boundary possibly non-smooth) promises to be utterly nebulous indeed, see also the concluding Remark below for some additional comments. In this paper, we shall instead start from the penalized approximate process (concerning the penalization coefficient $\beta_{0}$ see \eqref{Beta0} below)
\begin{equation}\label{AproxDif} 
X^{x}_{n}(t)=x+\int_{0}^{t}b(X^{x}_{n}(s))ds+\int_{0}^{t}
\sigma(X^{x}_{n}(s))dB(s)-n\int_{0}^{t}\beta_{0}(X^{x}_{n}(s))ds,
\end{equation}
proceed with a bona fide pathwise directional derivative and show, thanks to a tightness argument on $(X^{x}_{n},J^{x,\nu}_{n})$, where $J^{x,\nu}_{n}$ stands for the limit as $\epsilon\rightarrow 0$ in the ratio in \eqref{Ratio} but with $X^{x}_{n}$ instead of $X^{x}$, that the limit processes $(X^{x},J^{x,\nu})$ can be used at least to give a probabilistic representation for the gradient $\partial u(t,x)$ in \eqref{EqParbol}, even when $D$ is not convex by localization. 

Thus, the inability to deal with the full limit in \eqref{Ratio} actually shifts into an \textit{advantage}: we \textit{do not} need it to compute gradients!. Even when the full limit \eqref{Ratio} can be carried out, as in \cite{And2} for example, the formula \eqref{Ratio} for the gradient is costly for computer simulations as can be seen from Theorem 2.5 therein. Concrete simulations seem directly implementable: consider, for some $\epsilon>0$, the finite set of excursions (of $X$) $e_{\epsilon}(n)$, $n$ an integer, of duration $l(e_{\epsilon}(n))>\epsilon$. The Jacobian, from the start of an excursion, evolves continuously as in \eqref{JacobianFree} up to its collapse at the right end. Our work also provides some insights concerning \cite{Air}, \cite{IkeWat} and \cite{Pil}, see also additional comments in the concluding Remarks. Overall, our results and techniques are very different from the works listed above. The regularity of our data seems a little strong, which compensates however with the novelty of our results, and will be relaxed in a forthcoming paper.
\subsection{Notation}
We follow \cite{RevYor} for general probability notations, e.g. c\`{a}dl\`{a}g means right continuous with left-hand limit. $\mathcal{C}^{k}_{b}(\cdot)$ is the set of k times differentiable functions with bounded derivatives. We adopt the (Einstein) summation convention over repeated indices. $TD$ is the tangent bundle over the manifold $D$. $e_{i}$, $i=1,...,d$, is the canonical basis of $\mathbb{R}^{d}$. Inner product of $x$ and $y$ is $(x,y)$ or simply $x\cdot y$; for a quantity $c$, $c x\cdot y$ is the product $c$ times $x\cdot y$. For $\alpha\in \partial D$, $\gamma(\alpha)$ is the unit inward normal at $\alpha$ (regarded as a column vector $\gamma^{i}$), $N(\alpha)$ is projection onto the span of $\gamma(\alpha)$ and $N^{\bot}(\alpha)$ is the projection onto the tangent space $T_{\alpha}\partial D$. For a scalar function $f$, the (flat Euclidean) gradient is $\partial f$ or $\nabla f$ (the latter is preferred in a geometric setting); the Hessian matrix is $\partial^{2} f=[\partial^{2}_{x^{i}x^{j}} f]$. For a vector function $f=(f^{1},...,f^{d})$ and a vector $\eta$, $ (\partial f^{\cdot},\eta)$ is the vector with components $((\partial f^{1},\eta),(\partial f^{2},\eta),...)$. Unimportant constants $c,c^{\prime},...$ may change values while proofs are in process.

\subsection{Main result}
In addition to our assumptions on the coefficients $( b,\:\sigma)$ (see the beginning of this introductory section just below equation \eqref{Xdiffus}), we assume $D$ to be smooth, bounded, convex and uniformly non characteristic for $\sigma$, i.e. for some constant $c_{a}>0$, in some tubular neighborhood of $\partial D$, we have $\sum_{j}(\gamma(x)\cdot \sigma_{j}^{  \cdot}(x))^{2}\geq c_{a}$. The gradient of the squared distance function
\begin{equation}\label{Beta0} 
\beta_{0}(x)=(1/2)\partial(d(x,D)^{2})
\end{equation}
is used for penalization, see \cite{RenWu}, and the time spent outside $D$ by the process $X^{x}_{n}$ up to $t$
\begin{equation}\label{Tn} 
\mathcal{T}^{x,n}(t)=\int_{0}^{t}\mathbb{I}_{\lbrace X_{n}^{x}(s)\in D^{c}\rbrace}ds,
\end{equation}
is shown in \cite{BofBen} to satisfy 
\begin{equation}\label{BfBn}
\mathbb{E}\left(\mathcal{T}^{x,n}(t)\right)^{4}\leq  cn^{-2}    
\end{equation}
for some constant $c=c(t)$ independent of $x,\: n$ . Although $\beta_{0}(\cdot)$ is only Lipschitz ($ b,\:\sigma$ are in $\mathcal{C}^{1}_{b})$, we can still make first order derivatives in the initial condition $x$, by \cite{BouHir}, and our process below is a representative of the derivative process. Moreover, by \cite{AmbMan} on further derivatives on $\beta_{0}(\cdot)$, we can write for any direction $\nu$
\begin{equation}\label{JacobiAprox}
\begin{split}
&J^{x,\nu}_{n}(t)=\nu+\int_{0}^{t}(\partial b^{\cdot}(X_{n}^{x}(s)),J^{x,\nu}_{n}(s))ds+\int_{0}^{t}(\partial \sigma^{\cdot}_{ j}(X_{n}^{x}(s)),J^{x,\nu}_{n}(s))dB^{j}(s)\\
&-n\int_{0}^{t}N(X_{n}^{x}(s))J^{x,\nu}_{n}(s)d\mathcal{T}^{x,n}(s)
\end{split} 
\end{equation}
recall that $N(\alpha^{\prime})$, $\alpha^{\prime}\in D^{c}$, is the orthogonal projection on (the span of) the normal direction at $\alpha$ the projection of $\alpha^{\prime}$ on $\partial D$. Letting $e_{n}^{-}$ stand for an excursion of the process $X_{n}^{x}(t)$ in $D^{c}$ of duration $l(e_{n}^{-})>1/n$, we can remove the singularity $n$ and write the fundamental relation (paying attention to measurability matters)
\begin{equation}\label{JacobiNearFin}
\begin{split}
&J^{x,\nu}_{n}(t)\simeq  \nu+\int_{0}^{t}(\partial b^{\cdot}(X_{n}^{x}(s)),J^{x,\nu}_{n}(s))ds+\int_{0}^{t}(\partial \sigma^{\cdot}_{ j}(X_{n}^{x}(s)),J^{x,\nu}_{n}(s))dB^{j}(s)\\
&-\sum_{e_{n}^{-}}(NJ_{n}^{x,\nu})(t(e_{n}^{-}))
\end{split}
\end{equation}
where the sum is over those $e_{n}^{-}$ that are accomplished before $t$ and $t(e_{n})$ is some random time within these excursions. We now proceed with a rigorous derivation of this relation. Let us observe that any tightness criterion in the standard Skorohod metric topologies on the space $\mathbb{D}[0,t]$ of $\mathbb{R}^{d}$-valued c\`{a}dl\`{a}g functions on $[0,t]$ will not work since the limit points will ultimately be continuous, see \cite{Bil}. It turns out that the Jakubowski S-topology suits our needs, see \cite{Jak1} and \cite{Jak2}. We postulate a limit  
process that solves a key martingale problem.

The discontinuity structure of the limit process $J^{x,\nu}(t)$ is monitored by the Hsu \cite{Hsu1} point process of excursions of $X$ (from the boundary) which we now describe. Starting from $x\in D$, the law of the time and space exit from $D$ is absolutely continuous with density $(1/2)\partial _{\gamma(\alpha
)}p^{D}\left( t,\alpha ,x\right)$ ($p^{D}$ is the transition density of the killed process $X^{x}$ at $\partial D$). Let $
e=e(\alpha ,\beta )$ be the excursion which starts
at $\alpha $ and ends at $\beta $, and set $W^{\alpha }$ for the $e$'s that
start from $\alpha\in \partial D$. The point process of excursions, say $\mathcal{N}$, is naturally ordered in the local time scales. Let $S(\tau)$ be the right continuous inverse of the local time, i.e. for all $\tau\geq 0$
\begin{equation}\label{S}
S(\tau)=inf\lbrace t|L(t)>\tau\rbrace;    
\end{equation}
the trace $X^{x}(S(\tau))$ of $X^{x}$ on $\partial D$ is noted $\hat{X}(\tau)$; to each jump time $\tau >0$ of $S$ we have the excursion $e_{\tau }(t)=X(S(\tau -)+t)$ if $t\leq l(e_{\tau })$\ and $e_{\tau }(t)=X(S(\tau ))$ if $t\geq l(e_{\tau })$. $\mathcal{N}$ is quasi left-continuous and its compensating measure is explicitly given by (see \cite{Hsu1} p. 251)
\begin{equation*}
\widehat{\mathcal{N}}((0,\tau ]\times E)=\int_{0}^{\tau }Q^{\hat{X}(\theta) }\left(E\cap \lbrace  e(0)=\hat{X}(\theta)\rbrace \right)
d\theta ,
\end{equation*}
where $E$ is a measurable subset of $W=\cup _{\alpha }W^{\alpha }$ and the
excursion law $Q^{\alpha }$ on $W^{\alpha }$\ is given by
\begin{equation*}
Q^{\alpha }(e(t)\in dx,l(e)>t)=\frac{1}{2}\partial _{\gamma(\alpha
)}p^{D}\left( t,\alpha ,x\right) dx.
\end{equation*}
Here is our main result.
\begin{theorem}\label{ThmMain} 
Given our assumptions on the coefficients and $D$, let $L$ be the generator $L=L^{X}+L^{J}$, where
\begin{equation}\label{GenrtrCoupl} 
\begin{aligned}
& L^{X}{=}\dfrac{1}{2}a^{ij}(x)\partial _{x^{i}x^{j}}^{2}+ b^{i}(x)\partial_{x^{i}}, \\
& L^{J} {=}\dfrac{1}{2}(\sum_{k=1}^{d}( \partial\sigma_{k}^{i}(x),\nu)( \partial\sigma_{k}^{j}(x),\nu))\partial_{\nu^{i}\nu^{j}}^{2}\\
&+\sum_{k=1}^{d}\sigma_{k}^{i}(x)(\partial\sigma_{k}^{j}(x),\nu))\partial_{x^{i}\nu^{j}}^{2}+(\partial b^{i}(x),\nu)\partial_{\nu^{i}}
\end{aligned}
\end{equation}
with domain $\mathcal{D}$ the set of functions $F(x,\nu)\in \mathcal{C}^{2}_{b}(D\times\mathbb{R}^{d})$ subject to the following boundary conditions: for all $\nu\in\mathbb{R}^{d}$ we have \begin{equation}\label{BdryCondF1}
\partial_{\gamma}F(\cdot,\nu)(\alpha)=0, \:\forall \alpha\in \partial D,
\end{equation}
and  
\begin{equation}\label{BdryCondF2}
F(\alpha,\nu)=F(\alpha,N^{\bot}(\alpha)\nu), \:\forall (\alpha, \nu)\in \partial D\times\mathbb{R}^{d}.
\end{equation} 
Let $x\in D$ and $\nu$ be a (unit) direction in $\mathbb{R}^{d}$. The family of Markov processes $(X^{x}_{n},J^{x,\nu}_{n})$, $n$ an integer, where the first component is the penalization approximation of $X$ in \eqref{AproxDif} and the second component is the Jacobian in \eqref{JacobiAprox}, is tight in the Jakubowski S-topology and its limit points $(X^{x},J^{x,\nu})$ solve the martingale problem for $L$ in the sense that for all $F\in \mathcal{D}$
\begin{equation}
F(X^{x}
(t),J^{x,\nu}(t))-F(x,\nu)-\int_{0}^{t}LF(X^{x}(s),J^{x,\nu}(s))ds 
\end{equation}
is a $\mathbb{P}^{x,\nu}$-martingale. When this martingale problem is uniquely solvable, the limit $J^{x,\nu}$ is the c\`{a}dl\`{a}g semimartingale with a countable number of jumps which satisfies the SDE
\begin{equation}\label{EqJ}
\begin{aligned}
& J^{x,\nu}_{s}=\nu +\int_{0}^{t}(\partial b^{\cdot}(X_{s}^{x}),J^{x,\nu}_{s})ds+\int_{0}^{t}(\partial \sigma^{\cdot}_{ j}(X^{x}_{s})\cdot J^{x,\nu}_{s})dB^{j}_{s}\\
& -\int_{0}^{L(t)}\int_{W}
J^{x,\nu}(S(\tau)-)\cdot \gamma ( X^{x}(S(\tau)))
\gamma( X^{x}(S(\tau)))\mathcal{N}(d\tau,de).
\end{aligned}
\end{equation}
The matrix $J^{x}(t)$ with columns $J^{x,\nu=e_{j}}(t)$ (and $J^{x}(0)=\mathbb{I}$) where $\mathbb{I}$ is the unit $d\times d$ matrix, is an operator multiplicative functional in the sense that for all $0\leq r\leq s\leq t$
\begin{equation}\label{MOF}
J^{x,r}_{s}J^{x,s}_{t}=J^{x,r}_{t}.
\end{equation}
\end{theorem}

\begin{remark}
As soon as $X^{x}$ reaches $\partial D$, $J^{x,\nu}$ undergoes a discontinuity and restarts at $\zeta$ with $N^{\bot}(X^{x}(\zeta))H(\zeta)$ where for all $t>0$
\begin{equation}\label{H}
H(t)=\nu +\int_{ 0}^{t}\partial b^{\cdot}(X^{x}(s))\cdot J^{x,\nu}(s)ds+\int_{ 0}^{t}\partial \sigma^{\cdot}_{ j}(X^{x}(s))\cdot J^{x,\nu}(s)dB^{j}(s).
\end{equation}
The set of discontinuities of $J^{x,\nu}$, in real time, is the set of the right ends of the excursions intervals, which are stopping times. In local time scales, these are the times of the form $\lbrace S(\tau)|\tau\in \mathcal{J}(t)\rbrace$ where
\begin{equation}\label{Jumps}
\mathcal{J}(t)=\lbrace \tau|S(\tau-)\neq S(\tau)\rbrace.   
\end{equation}
\end{remark}

\section{Proof of theorem \ref{ThmMain}}
Note that thanks to the It\^{o} formula and to a standard Gronwall-Bellman argument, the solution of equation \eqref{EqJ} is unique for all $\nu$. Let $r\leq t$, by \cite{Die}, with the obvious notation, the matrix solution of
\begin{eqnarray*}
& J^{x,r}_{t}=\mathbb{I}+\int_{r}^{t}\partial b(X^{x}_{s})J^{x,r}_{s}ds+\int_{r}^{t}\partial \sigma_{ j}(X^{x}_{s}) J^{x,r}_{s}dB^{j}_{s}\\
&-\int_{0}^{L(t)}\int_{W}N(X^{x}_{S(\tau)})J^{x,r}_{S(\tau)-}\mathcal{N}(d\tau,de)    
\end{eqnarray*}
and is an operator multiplicative functional. 
\subsection{Tightness}\label{Tightness} 
The superscript $x$ is sometimes dropped when no ambiguity arises. It clearly suffices to deal with the component $J^{\nu}_{n}$. As was observed above, we need tightness in the Jakubowski S-topology for which it suffices to check that
\begin{equation}\label{CondVar} 
\sup_{n}\left(\sup_{s\leq t} \mathbb{E}\Vert J^{\nu}_{n}(s)\Vert+CV(t) \right)<\infty
\end{equation}
where $CV(\cdot)$ is the conditional variation, see \cite{MeyZhe} or Lemma \ref{CVBound} below. We shall need the fourth moment of the norm of $J^{\nu}_{n}$ and the following form of the Gronwall-Bellman inequality.
\begin{lemma}\label{Gronw}
Let $f$ be a continuous function on $[0,\infty)$ and $c_{1},c_{2}$ positive constants, if
\begin{equation*}
 f(t)\leq c_{1}+c_{2}(1+t)\int_{0}^{t}f(s)ds,   
\end{equation*}
then there are constants $c,c^{\prime},c^{\prime\prime}$ s.t.
\begin{equation*}
 f(t)\leq c\exp (c^{\prime}t^{2}+c^{\prime\prime}t).
 \end{equation*}
\end{lemma}
\begin{proof}
 Since $d/dt(c_{1}+c_{2}(1+t)\int_{0}^{t}f(s)ds)=c_{2}\int_{0}^{t}f(s)ds+c_{2}(1+t)f(t)$, we have
\begin{eqnarray*}
&& \frac{c_{2}\int_{0}^{t}f(s)ds+c_{2}(1+t)f(t)}{c_{1}+c_{2}(1+t)\int_{0}^{t}f(s)ds}\leq c_{2}(1+t) + \frac{c_{2}\int_{0}^{t}f(s)ds}{c_{1}+c_{2}(1+t)\int_{0}^{t}f(s)ds}\\
&& \leq c_{2}(1+t) +\frac{1}{1+t}  
\end{eqnarray*}
whence the result since for all $t$ we have $\log(1+t)\leq t$.
\end{proof}
\begin{lemma}\label{MomentEven} 
Given our assumptions, for all $t>0$ there exists a $c(t)$ s.t. $\forall n$ 
\begin{equation}\label{L2BoundOnJ} 
\mathbb{E}(\underset{0\leq s\leq t}{sup}\Vert J_{n}^{\nu}(s)\Vert^{4})\leq c(t).
\end{equation}
\end{lemma}
\begin{proof}
We have by It\^{o}
\begin{eqnarray*}
&&\Vert J_{n}^{\nu}(t)\Vert^{2}+2n\int_{0}^{t}\Vert N(X_{n}(s))J_{n}^{\nu}(s)\Vert^{2}d\mathcal{T}^{n}(s)=1\\
&& +2\int_{0}^{t}(\partial b^{\cdot}(X_{n}(s))\cdot J_{n}^{\nu}(s),J_{n}^{\nu}(s))ds+2\int_{0}^{t}(\partial\sigma^{\cdot}_{j}(X_{n}(s))\cdot J_{n}^{\nu}(s),J_{n}^{\nu}(s))dB^{j}(s)\\
&& +\sum_{j}\int_{0}^{t} \Vert \partial \sigma^{\cdot}_{j}(X_{n}(s))\cdot J_{n}^{\nu}(s)\Vert^{2}ds.
\end{eqnarray*}
Hence, by our conditions on the coefficients
\begin{equation*}
\Vert J_{n}^{\nu}(t)\Vert^{4}\leq c +ct\int_{0}^{t}\Vert J_{n}^{\nu}(s)\Vert^{4}ds+c \Vert\int_{0}^{t}(\partial\sigma^{\cdot}_{j}(X_{n}(s))\cdot J_{n}^{\nu}(s),J_{n}^{\nu}(s))dB^{j}(s)\Vert^{2}.   \end{equation*}
From now on, by Burkholder-Davis-Gundy and Lemma \ref{Gronw}, the rest of the proof is standard (it is not difficult to justify taking the expected value on the stochastic integral).
\end{proof} 
Before we pass on to the conditional variation, we need the following separation of scales key result. Set $S^{n}(\tau)$, $\tau\geq 0$, for the right continuous inverse of $\mathcal{T}^{n}(t)$, see \eqref{S}, $t\geq \zeta$ (recall $\zeta$ is the hitting time of $\partial D$ and is finite a.s.); $S^{n}(\tau)$ is allowed to make a jump at zero and start from $\zeta$. These are stopping times relative to $\mathcal{F}_{S^{n}(\tau)}$. Recall that $N(\alpha)$ refers to orthogonal projection on the normal field at $\alpha\in\partial D$, since $D$ is convex $N(\cdot)$ extends to the whole (half)-line spanned by $\gamma(\alpha)$.
\begin{lemma}\label{EdoForJ_N} 
For all $\tau\geq 0$, set $\hat{X}_{n}(\tau)=X_{n}(S^{n}(\tau))$ and $\hat{J}^{\nu}_{n}(\tau)=J^{\nu}_{n}(S^{n}(\tau))$. We have for all $\tau$
\begin{equation*}
N(\hat{X}_{n}(\tau))\hat{J}^{\nu}_{n}(\tau)=(N(\hat{X}_{n}(0))\hat{J}^{\nu}_{n}(0) )exp-n\tau+\int_{0}^{\tau}exp(-n(\tau-\theta))\hat{H}^{n}(\theta)d\theta
\end{equation*}
where $\hat{H}^{n}(\theta)=-N^{\bot}(\hat{X}_{n}(\theta))\hat{J}^{\nu}_{n}(\theta)+H^{n}(S^{n}(\theta))$ in which $H^{n}(t)$ is defined as in \eqref{H} with ${X}_{n}$ (resp. ${J}^{\nu}_{n}$) instead of ${X}$ (resp. instead of ${J}^{\nu}$). 
\end{lemma}
\begin{proof}
From \eqref{JacobiAprox} we clearly have the c\`{a}dl\`{a}g $\mathbb{R}^{d}$-valued random ODE
\begin{equation*}
N(\hat{X}_{n}(\tau))\hat{J}^{\nu}_{n}(\tau)=-N^{\bot}(\hat{X}_{n}(\tau))\hat{J}^{\nu}_{n}(\tau)+H^{n}(S^{n}(\tau))-n\int_{0}^{\tau}N(\hat{X}_{n}(\theta))\hat{J}_{n}^{\nu}(\theta)d\theta.
\end{equation*}
Hence, for continuity points $\tau$ (which correspond to the open downward excursions) the result is elementary in ODE theory. When $\tau$ is a jump time, the result still holds by right continuity.
\end{proof}
As a result of this, we have the
\begin{lemma}\label{CVBound} 
The sequence of processes $J_{n}^{\nu}(s)$, $0\leq s\leq t$ and $n$ an integer, is tight in the (Meyer-Zheng)-Jakubowski topology.
\end{lemma}
\begin{proof}
Let us first deal with the singular term. Using an integration by parts we have by Lemma \ref{EdoForJ_N} (with the obvious notation)
\begin{equation}\label{SingulrTerm}
\begin{aligned}
& n\int_{0}^{\mathcal{T}^{n}(t)}N(\hat{X}_{n}(\tau))\hat{J}_{n}^{\nu}(\tau) d\tau=n\lbrace \zeta<t\rbrace\int_{0}^{\mathcal{T}^{n}(t)}N(\hat{X}_{n}(\tau))\hat{J}_{n}^{\nu}(\tau) d\tau\\
& =\lbrace \zeta<t\rbrace N(X_{n}(\zeta))H^{n}(\zeta)(1-exp-n\mathcal{T}^{n}(t))\\
& +n\lbrace \zeta<t\rbrace\mathcal{T}^{n}(t)\int_{0}^{\mathcal{T}^{n}(t)}\hat{H}^{n}(\theta)exp-n(\mathcal{T}^{n}(t)-\theta)d\theta\\
&-n\lbrace \zeta<t\rbrace\int_{0}^{\mathcal{T}_{n}(t)}\tau\hat{H}^{n}(\tau)d\tau + n^{2}\lbrace \zeta<t\rbrace\int_{0}^{\mathcal{T}_{n}(t)}(\tau\int_{0}^{\tau}\hat{H}^{n}(\theta)exp-n(\tau-\theta)d\theta)d\tau.
\end{aligned}
\end{equation}
We have
\begin{equation*}
\lbrace \zeta<t\rbrace N(X_{n}(\zeta))H^{n}(\zeta)=\lbrace \zeta<t\rbrace N(X_{n}(\zeta))H^{n}(\zeta\wedge t),
\end{equation*}
so that by the boundedness of the coefficients, Jensen's inequality and Lemma \ref{MomentEven}
\begin{equation*}
\mathbb{E}\Vert\lbrace \zeta<t\rbrace N(X_{n}(\zeta))H^{n}(\zeta)\Vert\leq \mathbb{E}\Vert 
H^{n}(\zeta\wedge t)\Vert\leq c. 
\end{equation*}
Next, in a similar way, by the relation \eqref{BfBn} we have 
\begin{eqnarray*}
&& n\mathbb{E}\lbrace \zeta<t\rbrace\int_{0}^{\mathcal{T}_{n}(t)}\tau\Vert \hat{H}^{n}(\tau)\Vert d\tau\leq cn\mathbb{E}\lbrace \zeta<t\rbrace\underset{\tau\leq \mathcal{T}^{n}(t)}{sup} \Vert \hat{H}^{n}(\tau)\Vert (\mathcal{T}^{n}(t))^{2}\\
&& \leq cn\mathbb{E}\lbrace \zeta<t\rbrace\underset{\zeta\leq s\leq t}{sup} \left(\Vert J^{\nu}_{n}(s) \Vert +\Vert               H^{n}(s)\Vert \right)(\mathcal{T}^{n}(t))^{2}\\
&& \leq c\mathbb{E}(\underset{s\leq t}{sup} \Vert J_{n}^{\nu}(s)\Vert^{2})+ c\mathbb{E}(n^{2}\mathcal{T}^{n}(t)^{4})\leq c<\infty, 
\end{eqnarray*}
uniformly in $n$. 

The second term on the right hand-side in equation \eqref{SingulrTerm} is easy (under an expectation, it is either an asymptotically small term or just bounded by a constant) and the last term is treated using the elementary inequality $n\int_{0}^{\tau}exp-n(\tau-\theta)d\theta\leq 2$. To sum up
\begin{equation*}
n\mathbb{E}\int_{0}^{\mathcal{T}_{n}(t)}\Vert N(\hat{X}_{n}(\tau))\hat{J}_{n}^{\nu}(\tau)\Vert d\tau\leq c,
\end{equation*}
uniformly in $n$. Now let $t_{k}$, $0=t_{0}<t_{1}<...<t_{m}<t$, be a subdivision of $[0,t]$. We have uniformly in $n$
\begin{eqnarray*}
&&CV(t)=\sum_{k} \mathbb{E}\Vert \mathbb{E}\left(J_{n}^{\nu}(t_{k+1})-J_{n}^{\nu}(t_{k})|\mathcal{F}_{t_{k}}\right)\Vert \leq \int_{0}^{t}\mathbb{E}\Vert(\partial b^{\cdot}(X_{n}(s)),J_{n}^{\nu}(s))\Vert ds\\
&&+n\mathbb{E}\int_{0}^{t}\Vert N(X_{n}(s))J_{n}^{\nu}(s)\Vert d\mathcal{T}^{n}(s)\leq c.
\end{eqnarray*}
\end{proof}
\subsection{The limit points}
We use localization. Let us first take $D$ to be the upper half-space $\mathbb{R}^{d-1}\times\mathbb{R}^{+}$ and Let $F\in\mathcal{D}$. Consider in $D$ the process $\tilde{X}_{n}$ with horizontal component $(X_{n}^{1},...,X_{n}^{d-1})$ and vertical one $\tilde{X}^{d}_{n}=X^{d,+}_{n}$. Note that $X_{n}$ differs from $\tilde{X}_{n}$ only during the small downward excursions. By the Tanaka formula
\begin{equation*}
\tilde{X}^{d}_{n}(t)=x^{d}+\int_{0}^{t}\mathbb{I}_{\lbrace X_{n}^{d}(s)>0\rbrace}b^{d}(X_{n}(s))ds+\int_{0}^{t}\mathbb{I}_{\lbrace X_{n}^{d}(s)>0\rbrace}\sigma^{d}_{j}(X_{n}(s))dB^{j}(s)+\dfrac{1}{2}L^{n}(t),
\end{equation*}
where $L^{n}$ is the (semimartingale) local time of $X^{d}_{n}$ at zero. Set $Y^{n}=(\tilde{X}_{n},J_{n}^{\nu})$; by the boundary conditions \eqref{BdryCondF1} and \eqref{BdryCondF2} we have in particular $ \partial_{\gamma}F(\alpha,\cdot)(\nu)=0$ and by the It\^{o} formula we have for $r\leq t$  
\begin{eqnarray*}
&&F(Y^{n}_{t})-F(Y^{n}_{r})\\
&&=\int_{r}^{t}\partial_{x} F(Y^{n}_{s}) \cdot b(X_{n}(s))ds+\int_{r}^{t} ( \partial_{\nu} F(Y^{n}_{s}) ,\partial b^{\cdot}(X_{n}(s))\cdot J_{n}^{\nu}(s) ) ds\\
&&+\int_{r}^{t}\partial_{x} F(Y^{n}_{s})\cdot \sigma_{j}^{\cdot}(X_{n}(s)) dB^{j}_{s} +\int_{r}^{t}(\partial_{\nu} F(Y^{n}_{s}),\partial \sigma_{j}^{\cdot}(X_{n}(s))\cdot J_{n}^{\nu}(s)) dB^{j}_{s}\\
&&+\dfrac{1}{2}\int_{r}^{t}\partial_{x^{k}x^{l}}^{2}F(Y^{n}_{s})a^{kl}(X_{n}(s)) ds\\
&&+\int_{r}^{t}\partial_{x^{k}\nu^{l}}^{2}F(Y^{n}_{s})(\sum_{j=1}^{d}\sigma_{j}^{k}(X_{n}(s))\partial\sigma_{j}^{l}(X_{n}(s))\cdot J_{n}^{\nu}(s))ds\\
&&+\dfrac{1}{2}\int_{r}^{t}\partial_{\nu^{k}\nu^{l}}^{2}F(Y^{n}_{s})(\sum_{j=1}^{d}(\partial\sigma_{j}^{k}(X_{n}(s)), J_{n}^{\nu} (s))( \partial\sigma_{j}^{l}(X_{n}(s)),J_{n}^{\nu} (s)))ds\\
&&+R^{n}(r,t)
\end{eqnarray*}
where
\begin{eqnarray*}
&& R^{n}(r,t)=-\int_{r}^{t}\mathbb{I}_{\lbrace X_{n}^{d}(s)\leq 0\rbrace}\partial_{x^{d}} F(Y^{n}_s) (b^{d}(X_{n}(s))ds+\sigma_{j}^{d}(X_{n}(s)) dB^{j}_{s})\\
&& -\int_{r}^{t}\mathbb{I}_{\lbrace X_{n}^{d}(s)\leq 0\rbrace}(\sum_{k<d}\partial_{x^{k}x^{d}}^{2}F(Y^{n}_{s})a^{kd}(X_{n}(s)) )ds\\
&&-\dfrac{1}{2}\int_{r}^{t}\mathbb{I}_{\lbrace X_{n}^{d}(s)\leq 0\rbrace}\partial_{x^{d}}^{2}F(Y^{n}_{s})a^{dd}(X_{n}(s))ds\\
&& -\int_{r}^{t}\mathbb{I}_{\lbrace X_{n}^{d}(s)\leq 0\rbrace} \left(\sum_{l=1}^{d} \partial_{x^{d}\nu^{l}}^{2}F(Y^{n}_{s})(\sum_{j=1}^{d}\sigma_{j}^{d}(X_{n}(s))\partial\sigma_{j}^{l}(X_{n}(s))\cdot J_{n}^{\nu}(s))\right)ds.
\end{eqnarray*}
Let $L^{X}(x)$ and $L^{J}(x,\nu)$ be the generators from equation \eqref{GenrtrCoupl} and let $\varphi_{r}(\cdot)$ be a  abounded continuous functional over the Skorohod space $\mathbb{D}[0,t]$ which depends only on the past up to $r$, we have 
\begin{eqnarray*}
&& \mathbb{E}[(F(Y_{n}(t))-F(Y_{n}(r))
-\int_{r}^{t}(L^{X}+L^{J})(X_{n}(s),J_{n}^{\nu}(s))
F(Y_{n}(s))ds)\varphi_{r}(\cdot)]\\
&& =\mathbb{E}(R^{n}(r,t)\varphi_{r}(\cdot)).
\end{eqnarray*}
By our conditions on $F(x,\nu)$ and on our coefficients, we clearly have by well known inequalities and Lemma \ref{MomentEven} 
\begin{equation*}
\mathbb{E}|R^{n}(r,t))|^{2}\leq c(t) \mathbb{E}(\mathcal{T}^{n}(t)+(\mathcal{T}^{n}(t))^{1/2}+(\mathcal{T}^{n}(t))^{1/4})\leq c(t) n^{-1/8}
\end{equation*}
and since
\begin{equation*}
\mathbb{E}\underset{s\leq t}{sup} \Vert \tilde{X}_{n}(s)-X(s)\Vert^{2}\leq c\mathbb{E}(\underset{s\leq t}{sup} \Vert X_{n}(s)-X(s)\Vert^{2})+c\mathbb{E}\underset{s\leq t}{sup} (X_{n}^{d}(s)-X^{d,+}_{n}(s))^{2},
\end{equation*}
then by Proposition 3.4 and Theorem 3.6 of \cite{RenWu} it now suffices to let $n\rightarrow\infty$, perhaps through a subsequence.

\section{A probabilistic formula for the gradient}\label{Apl}
As explained above, some geometric framework is needed. The domain $D$ is a flat Riemann manifold with boundary and we view the gradient vector $v(x)=\partial u(t,x)$ as rather a differential 1-form $v(x)=\partial_{x^{i}} u(t,x)dx^{i}$, i.e. an element of the cotangent bundle $T_{x}^{*}D$ (in what follows $d(\cdot)$ is the exterior derivative on differential forms); consider the Stratonovitch reflected diffusion
\begin{equation}\label{XStratono} 
 X^{x}_{t}=x+\int_{0}^{t}V_{j}(X^{x}_{s})\circ dB^{j}(s)+\int_{0}^{t}V_{0}(X^{x}_{s})ds
+\int_{0}^{t}\gamma(X^{x}_{s})dL^{x}(s),
\end{equation} 
where the coefficients $V_{j}=\sigma_{j}^{i}\partial_{i}$, $j=0,...,d$, are smooth vector fields. It has the It\^{o} form \eqref{Xdiffus} with $b^{i}=V^{i}_{0}+(1/2)\sum_{j=1}^{d}V_{j}(\sigma_{j}^{i})$. The generator of $X$ is the (suitably closed) sum of squares operator
\begin{equation}\label{SumSquar} 
\mathcal{L}(\cdot)=\sum_{j=1}^{d}\mathcal{L}_{V_{j}}^{2}(\cdot)+\mathcal{L}_{V_{0}}(\cdot)
\end{equation}
where $\mathcal{L}_{V_{j}}$ stands for the Lie derivative along the vector field $V_{j}$, $j=0,...,d$. Recall, see any textbook on differential geometry, that a Lie derivative is just ordinary directional derivative on a scalar function and that the Lie derivative, at $x$, of a $p$-form $w$ along a vector field $V$ is given by the limit of the ratio $(\phi_{t}^{\star}w_{\phi_{t}x}-w_{x})/t$ as $t\rightarrow 0$ where $\phi_{t}$ is the flow generated by the vector field $V$ and $\phi_{t}^{\star}$ is the associated pull-back operator on forms. Recall a form $w$ is closed if $\partial_{x^{j}}w_{i}=\partial_{x^{i}}w_{j}$ for all $i,j$ and exact if $w=d\varphi$ for some scalar function $\varphi$, moreover $d^{2}(\cdot)=0$.  

We shall need the Cartan formula for an exterior $p$-form
\begin{equation}\label{Cartan1} 
\mathcal{L}_{V}=i_{V}\circ d+d\circ i_{V}
\end{equation}
where $i_{V}$ is the interior product (i.e. a contraction and also, by the way, an anti-derivative) on a $p$-form $w$ given by $(i_{V}w)(U_{2},...,U_{p})=w(V,U_{2},...,U_{p})$ over a collection of $p-1$ arbitrary argument vector fields $(U_{2},...,U_{p})$. A consequence of equation \eqref{Cartan1} is that for a $1$-form $w$ we have the extremely useful formula (at $x$) 
\begin{equation}\label{Cartan2} 
dw(V_{x},V^{\prime}_{x})=V_{x}(w(V^{\prime}))+V^{\prime}_{x}(w(V))-w([V,V^{\prime}])
\end{equation}
where the bracket, i.e. the Lie derivative, $[V,V^{\prime}]$ is the vector field whose differential operator is the commutator (acting on scalar functions $f$) $V_{x}(V^{\prime}(f))-V^{\prime}_{x}(V(f))$. 

Given our assumptions, we know from Theorem 7 Chapter 9, Section 6, of \cite{Fri} that $u(t,\cdot)$ is at least in $\mathcal{C}^{3}$. Since exterior differentiation commutes with the Lie derivative, we see that the time dependent 1-form $v(t,x)=du(t,x)$ satisfies the equation
\begin{equation}\label{ParbolMathcalL}
\partial_{t}v(t,x)=\mathcal{L}v(t,x), \:x\in \overset{o}{D}    
\end{equation}
Regarding the boundary conditions, note that the Cauchy condition $u(0,x)=f(x)$ transforms into $v(0,x)=df_{x}$ and the Neumann condition in \eqref{EqParbol} implies the Dirichlet condition $v(t,\alpha)(\gamma)=0$, $\alpha\in \partial D$; therefore from the fact that $dv(t,x)=0$, we see (upon applying \eqref{Cartan2}) that $v$ satisfies the so called absolute boundary condition
\begin{equation}\label{AbsBdryCond}
w^{(N)}\overset{\partial D}{=}0, (dw)^{(N)}\overset{\partial D}{=}0, w \: being \: a \: 1-form
\end{equation}
where the superscript $N$ in a $p$-form stands for its "normal part" which is simply the same form but acting on the normal parts of their $p$ (contravariant) vector arguments at the boundary. Remark in passing that, in a general Riemann manifold, the whole mosaic of the different boundary conditions is masterly reviewed in \cite{MitTay}. For a probabilistic representation of the problem \eqref{ParbolMathcalL}, \eqref{AbsBdryCond} we shall proceed with a "forcing" and this requires help from Analysis rather than from Probability Theory since another derivative on $J^{x}(t)$ does not seem applicable. Therefore we assume that:

Condition A: The semigroup, for continuous $F(\cdot,\cdot)$, $\mathbb{E}^ {x,\nu}F(X(t),J(t))$ is at least $\mathcal{C}^{2}$ in $x$ for all $\nu$. 

At this stage, let us compute $\mathcal{L}$ in \eqref{SumSquar} for any closed $1$-form $w_{x}=w_{i}(x)dx^{i}$ in $\mathcal{C}^{2}$. The following result should be implicit in the literature about differential geometry. Since we couldn't find it explicitly written, we shall provide a proof. 
\begin{lemma}
Let $w$ be as above, and set $b^{i}=V^{i}_{0}+(1/2)\sum_{j=1}^{d}V_{j}(\sigma_{j}^{i})$, we have
\begin{equation}\label{GradSystm} 
\begin{aligned}
&\mathcal{L}w|_{k}{=}\dfrac{1}{2}a^{ij}(x)\partial _{x^{i}x^{j}}^{2}w_{k} + b^{i}(x)\partial_{x^{i}}w_{k}
+\dfrac{1}{2}(\partial_{x^{k}}a^{ij}(x))\partial _{x^{i}}w_{j}+(\partial_{x^{k}}b^{i}(x))w_{i}\\
& =L^{X}w_{k}+\dfrac{1}{2}(\partial_{x^{k}}a^{ij}(x))\partial _{x^{i}}w_{j}+(\partial_{x^{k}}b^{i}(x))w_{i}.
\end{aligned}
\end{equation}
\end{lemma}
\begin{proof}
Upon applying the relations \eqref{Cartan1} and \eqref{Cartan2}, it is not difficult to see (thanks to elementary derivations, cancellations and renaming of indices) that for all $j=0,...,d$
\begin{equation}\label{1stLie} 
\mathcal{L}_{V_{j}}w|_{k}=(\partial_{x^{k}}\sigma_{j}^{i})w_{i}+\sigma_{j}^{i}\partial_{x^{i}}w_{k}.
\end{equation}    
Applying again the above formula, thanks to similar calculations, it turns out (as expected) that the result amounts to taking formal derivatives on the system \eqref{EqParbol}. 
\end{proof}
We shall prove the following 
\begin{theorem}
Under the hypotheses of Theorem \ref{ThmMain}, suppose moreover that the condition A just above holds. Then the 1-form $v_{i}(t,x)=\mathbb{E}(\nabla f( X^{x}(t) )\cdot J^{x,e_{i}}(t))$ where $J^{x,e_{i}}(t)$ is defined in Theorem \ref{ThmMain} solves the system \eqref{ParbolMathcalL}, \eqref{AbsBdryCond} in the classical sense.
\end{theorem}
\begin{proof}
Set for the integer $n$
\begin{equation}\label{vn}
v^{n}(t)=\mathbb{E}(\nabla f(X^{x}_{n}(t))J_{n}^{x}(t))
\end{equation}
since $\Vert v_{i}^{n}(t)\Vert_{L^{2}}\leq c(t)$, $i=1,...,d$, by Lemma \ref{MomentEven}, up to an extraction, the sequence \eqref{vn} has a weak-star limit. It then suffices to apply \cite{Sco} to deduce that (by abuse of notation we keep writing $v$ as in the $v=du(t,x)$ where $u$ solves \eqref{EqParbol})
\begin{equation}\label{vEnfin}
 v(t,x)=\mathbb{E}(\nabla f(X^{x}(t))J^{x}(t))
 \end{equation}
is closed (indeed exact).

Now, let us first check that $v$ satisfies the PDE \eqref{GradSystm} inside $D$. We take a smooth $f\in \mathcal{D}^{X}$, the general case follows by a density argument. Let $x$ be inside $D$ and let $\zeta(\epsilon)$ be the first exit from a fixed small ball $B(x,\epsilon)$. By \eqref{MOF} and by the Markov property we have (for each component of $v$)
\begin{eqnarray*}
&&v(t+\Delta t,x)=\mathbb{E}^{x}[\nabla f(X_{t+\Delta t})J^{\Delta t}_{t+\Delta t}J^{0}_{\Delta t}]\\
&&=\mathbb{E}^{x}[\mathbb{E}^{X_{\Delta t}}(\nabla f(X_{t})J^{0}_{t})J^{0}_{\Delta t}].
\end{eqnarray*}
By the above It\^{o} formula, we have for all $\nu=e_{1},...,e_{d}$
\begin{eqnarray*}
&& \mathbb{E}^{x}( v(t,X_{\Delta t}),J_{\Delta t}^{\nu}) -( v(t,x),\nu)=\mathbb{E}^{x}\int_{0}^{\Delta t}( L^{X}(X_{s})v(t,X_{s}),J_{s}^{\nu}) ds\\
&&+\mathbb{E}^{x}\int_{0}^{\Delta t}v_i(t,X_{s})( \nabla b^i(X_{s}),J_{s}^{\nu}) ds\\
&&+\mathbb{E}^{x}\int_{0}^{\Delta t}\sum_{k}(\nabla v_i(t,X_{s}),\sigma^{\cdot}_{k}(X_{s}))( \nabla \sigma^{i}_{k}(X_{s}),J_{s}^{\nu}) ds\\
&& +\mathbb{E}^{x}\int_{0}^{\Delta t}\left( 
J_{s}^{\nu}, \nabla v_{\cdot}(t,X_{s}) \cdot \gamma(X_{s})\right) dL_{s}.
\end{eqnarray*}
The last term on the right hand-side above equals 
\begin{equation*}
\mathbb{E}^{x}\lbrace\zeta(\epsilon) \leq \Delta t\rbrace\int_{0}^{\Delta t}\left( 
J_{s}^{\nu}, (\nabla v_{\cdot}(t,X_{s}) \cdot \gamma(X_{s}))\right) dL_{s};   
\end{equation*}
once divided by $\Delta t$, it is not difficult to show that it does not contribute in the limit as $\Delta t\rightarrow 0$ by uniform ellipticity of $a(x)$, our assumptions on the coefficients and Lemma \ref{MomentEven}.

By right continuity we have 
\begin{eqnarray*}
&&(v(t+\Delta t,x)-v(t,x))/\Delta t\rightarrow\\
&&L^{X}(x)v_i(t,x)\delta^{i}_{\nu}+v_i(t,x)( \nabla b^i(x),\nu)\\
&&+\sum_{k}(\nabla v_i(t,x),\sigma_{k}^{\cdot}(x))( \nabla \sigma_{k}^{i}(x),\nu)\\
&&=L^{X}(x)v_{\nu}(t,x)+\sum_{i}\partial_{\nu}b^i(x)v_i(t,x)\\
&&+\sum_{ij}\partial_{x^j}v_i(t,x)(\sum_{k}\sigma_{k}^{j}(x)\partial_{\nu}\sigma_{k}^{i}(x)). 
\end{eqnarray*}
Now, we pass on to the boundary condition. As far as the Dirichlet condition is concerned, let $\alpha\in \partial D$. Let $g^{j}$, $j=1,...,d$, be an orthonormal basis for $T^{*}_{\alpha}D$ with $g^{1}$ given by $g^{1}(V)=(V,\gamma)$. Then $v^{(N)}=(v,g^{1})g^{1}$; hence $v^{(N)}=0$ if and only if $(v,g^{1})=0$.  That is we need only check $\mathbb{E}(\nabla f(X^{\alpha}(t)), J^{\alpha,\gamma}(t))=0$ for all $t$. Since we start from the boundary, the direction $\gamma$ in \eqref{EqJ} is immediately projected on $T_{\alpha}\partial D$ and we have by the It\^{o} formula
\begin{eqnarray*}
&&\Vert J_{t}^{\alpha,\gamma}\Vert^{2}=2\int_{0}^{t}(\nabla b^{\cdot}(X^{\alpha}_{s})\cdot J_{s}^{\alpha,\gamma},J_{s}^{\alpha,\gamma})ds+2\int_{0}^{t}(\nabla\sigma_{j}^{\cdot}(X^{\alpha}_{s})\cdot J_{s}^{\alpha,\gamma},J_{s}^{\alpha,\gamma})dB^{j}_{s}\\
&&+\underset{s\in \mathcal{J}(t)}{\sum}\left( \Vert J_{s}^{\alpha,\gamma}\Vert^{2}-\Vert J_{s-}^{\alpha,\gamma}\Vert^{2}\right).
\end{eqnarray*}
Since the last sum is clearly negative and since our coefficients are bounded, the Dirichlet condition follows immediately from a standard Gronwall-Bellman argument.

There remains to deal with the Neumann condition. Since the Dirichlet condition holds and since our form $v(t,x)$ in \eqref{vEnfin} is closed, the Neumann condition is automatic (by \eqref{Cartan2}). 
\end{proof}
\subsection{The Weingarten map}\label{Weingrtn}
Many works cited in the introductory Section, especially all those papers with a non constant diffusion matrix, do not consider a properly curved boundary. We now show how the Weingarten map, or shape operator, $\mathcal{S}(\cdot)$ on $\partial D$ appears naturally. This is the map $T_{\alpha}\partial D\ni \eta\mapsto \mathcal{S}(\alpha)(\eta)=-\nabla_{\eta}\gamma(\alpha)\in T_{\alpha}\partial D$. We have for $f\in \mathcal{D}(L^{X})$ 
\begin{equation*}
\begin{split}
&  ( \nabla\partial_{\alpha^{i}} f(\alpha),\gamma(\alpha))=\sum_{j}^{d}\left( 
\partial_{\alpha^{i}}(\partial_{\alpha^{j}}f(\alpha)\gamma^{j}(\alpha))-
\partial_{\alpha^{j}}f(\alpha) \partial_{\alpha^{i}}\gamma^{j}(\alpha)\right)\\
&=-\sum_{j}^{d}\partial_{\alpha^{j}}f(\alpha) \partial_{\alpha^{i}}\gamma^{j}(\alpha),
\end{split}
\end{equation*}
so that for fixed $\eta\in T_{\alpha}\partial D $
\begin{equation*}
\sum_{i}^{d}\eta^{i}( \nabla\partial_{\alpha^{i}} f(\alpha),\gamma(\alpha)) =-( \nabla f(\alpha),\partial_{\eta}n(\alpha)\rangle=-( \nabla f(\alpha),\mathcal{S}(\alpha)(\eta)).
\end{equation*}
Hence, by the It\^{o} formula the component $\nu=e_{i}$ of our gradient has a local time term 
\begin{equation}\label{WeingForf}
-\mathbb{E}\int_{0}^{t} (\nabla f(X_{s}),\mathcal{S}(X_{s})(J_{s}^{\nu})) dL(s)=-\mathbb{E}\int_{0}^{L(t)} (\nabla f(\hat{X}_{\tau}),\mathcal{S}(\hat{X}_{\tau})(\hat{J}^{\nu}_{\tau} ))d\tau,
\end{equation}
recall the Jacobian, for all $\nu$, is tangential at the times of the form $S(\tau)$, $\tau\geq 0$.
\section{Concluding remarks}
\begin{enumerate}
\item
Arguing as in \cite{IkeWat} p. 292, in boundary value problems in PDE theory (with \textit{arbitrary} smooth $D$) we often need to locally flatten the boundary in order to make certain estimates. Our method does provide one way to deal with such issues.
We are unable to exhibit a global Jacobian for a general diffusion when $D$ is not convex and the following picture is intuitively suggestive. Imagine, from the Jacobi field along $X^{x}$ point of view, a fluid of tiny sticks $\epsilon\nu$ (at time zero) sloshing inside our domain $D$. A stick whirls violently as one of its tips hurts prominent boundary bumps pointing inside the domain and the subsequent evolution is thereby quite intractable. Overall, it seems clear that there are less whirlings inside a convex body or inside a non convex one but with less diffusion as in the constant diffusion matrices of \cite{Air}, \cite{Bur} and \cite{And2}. Indeed, for a non convex body, let us first define an approximate Jacobian. Given some smoothness, a unique metric projection is defined within a finite shell around $\partial D$ in which penalization is smoothly carried out, see \cite{RenWu}. In order to make a further differentiation, the shell is also adjusted to fulfill Theorem 2.1 in \cite{AmbMan}. Namely, there is a uniform sufficiently small $r_{0}>0$ s.t. for all $x$ in the outer shell $D^{-}(3r_{0})=\lbrace x\in D^{c}: 0<d(x,D)\leq 3r_{0}\rbrace$, there is a unique boundary point noted $\pi(x)$, which is closest to $x$. Let us take a smooth $\rho(t)=t$ when $t\in[0, 4r_{0}^{2}]$, $\rho(t)=9r_{0}^{2}$ when $t\in[ 9r_{0}^{2},\infty)$ and monotone increasing in between. Now, formally take
\begin{equation}\label{Beta} 
\beta(x)=(1/2)\partial(\rho(d(x,D)^{2})).
\end{equation}
instead of $\beta_{0}$ in \eqref{AproxDif}. By the general theory we can write for a direction $\nu$
\begin{equation}\label{JacobiAprox+}
\begin{split}
&J_{n}^{\nu}(t)=\nu+R_{n}^{\nu}(t)+\int_{0}^{t}(\partial b^{\cdot}(X_{n}(s)),J_{n}^{\nu}(s))ds\\
&-n\int_{0}^{t}\frac{d\rho}{dt}(d^{2}(X_{n}(s)))N(X_{n}(s))J_{n}^{\nu}(s)d\mathcal{T}^{n}(s)
\end{split} 
\end{equation}
where 
\begin{equation}\label{Reste} 
R_{n}^{\nu}(t)=-2n\int_{0}^{t}\frac{d^{2}\rho}{dt^{2}}(d^{2}(X_{n}(s)))(\beta_{0}(X_{n}(s)),J^{\nu}_{n}(s))\beta_{0}(X_{n}(s))d\mathcal{T}^{n}(s).
\end{equation}
The proof of tightness, with no stochastic integrals on the right hand-side, runs as in Section \ref{Tightness} concerning the bounds on (even) moments of the approximate Jacobian except for a Bihari-Langenhop argument instead of a standard Gronwall-Bellman one due to the singular additional term $R^{\nu}_{n}$ which can't be sent to the left hand-side as in the proof of Lemma \ref{MomentEven}. We further need the $L^{p}$ moment estimates of \cite{RenWu} concerning the distance process $d(X_{n},\partial D)$ to deal with the singularity $n$. The details are left to the reader.
\item
Our result should also hold in a curved Riemann manifold. However, the right concept for drawing a diffusion on a manifold is not the Lie derivative but that of the covariant one, e.g. the Lie derivative lacks important linearity properties enjoyed by the covariant derivative. By now the ideas are clear but the details are not straightforward because, e.g., the Laplace-Beltrami operator can be written as the sum of squares \eqref{SumSquar} only locally; moreover we need to address the construction of the Jacobian, which should live on the tangent space at given $x$, and its related discontinuities (but these are of secondary difficulty since they occur on the flat tangent spaces). As indicated in \cite{ElwLi}, the Jacobian should be studied within the framework of the covariant SDEs and connections on the bundle $TD$ seem necessary.
\item
The trace Jacobian matrix process on $\partial D$ should be invertible and seems to provide a tool for studying fine properties of the boundary processes which are associated with pseudo-differential operators, especially in the non self-adjoint case. 
\item
Since our main ingredient \cite{BofBen} allows for some degeneracy, our method seems to be applicable to some extent to degenerate PDEs.
\item  
The last term in \eqref{EqJ} is the sum $-\sum_{\tau\in \mathcal{J}(t)}J^{\nu}(S(\tau)-)\cdot\gamma(X(S(\tau))) \gamma(X(S(\tau)))$, see \eqref{S}, and is essentially a random L\'{e}vy system over the extended boundary process $(S(\tau),X(S(\tau))$ and it is important to realize it does not depend on the excursions of $X$ inside $D$ but only on the boundary process. 

\end{enumerate}


$\,$

$\,$


\begin{thebibliography}{99}


\bibitem{Air}
Airault, H. (1976). ``Perturbations singuli\`{e}res et solutions stochastiques de probl\`{e}mes de D. Neumann-Spencer.'' \emph{J. Math. Pures Appl.} 55 (3): 233--267.

\bibitem{AmbMan}
Ambrozio, L., Mantegazza, C. (1998). ``Curvature and distance function from a manifold.'' \emph{J. Geom. Analysis} 8: 723--748. 

\bibitem{And1}
Andres, S. (2009). ``Pathwise differentiability for SDEs in a convex polyhedron with oblique reflection.'' \emph{Ann. I.H.P.} 45 (1): 104--116.

\bibitem{And2} 
Andres, S. (2011). ``Pathwise differentiability for SDEs in a smooth domain with reflection.'' \emph{Elec. J. Prob.} 28: 845--879.

\bibitem{ArnLi} 
Arnaudon, M., Li, X.~M. (2017). `` Reflected Brownian motion: selection, approximation and linearization.'' \emph{Elec. J. Prob.} 22: 1--55.

\bibitem{Bil} 
Billingsley, P. (1999). \emph{Convergence of probability measures.} Wiley. 

\bibitem{BofBen} 
Boufelgha, N., Bench\'{e}rif-Madani, A. (2023). ``A penalization limit theorem for the boundary local time of a reflected diffusion.'' \emph{Markov Proc. Rel. Fields} 22 (38): 41--61.

\bibitem{BosCisTal} 
Bossy, M., Ciss\'{e}, M., Talay, D. (2011). ``Stochastic representation of derivatives of solutions on one-dimensional parabolic variational inequalities with Neumann boundary conditions.'' \emph{Ann. I.H.P. Prob. Stat.} 47 (2): 395--424.

\bibitem{BouHir} 
Bouleau, N., Hirsch, F. (1989). ``On the derivability, with respect to initial data, of the solution of a stochastic differential equation with lipschitz coeffcients.'' \emph{S\'{e}minaire de Th\'{e}orie du Potentiel, Paris} 9: 39--57.

\bibitem{Bur} 
Burdzy, K. (2009). ``Differentiability of stochastic flow of reflected Brownian motions.'' \emph{Elec. J. Prob.} 14 (75): 2182--2240.

\bibitem{BurLee} 
Burdzy, K., Lee, J. ~M. (2010). ``Multiplicative functional for reflected Brownian motion via deterministic ODE.'' \emph{Illinois. J. Math.} 54 (3): 895--925. 

\bibitem{DoeZam} 
 Deuschel, J.~D., Zambotti, L. (2005). ``Bismut-Elworthy's formula and random walk representation for SDEs with reflection.'' \emph{Stoch. Process. Appl.} 115; 907--925.

\bibitem{Die} 
Dietz, H.~M. (1991). ``On the solution of matrix-valued linear SDEs driven by semimartingales.'' \emph{Stochastics and Stochastic reports} 34: 127--147.

\bibitem{ElwLi}
Elworthy, K.~D., Li, X.~M., (1994). ``Formulae for the Derivatives of Heat Semigroups.'' \emph{J. Func. Anal.} 125 (1):  252--286.  

\bibitem{Fri} 
Friedman, A. (1964). \emph{Partial differential equations of parabolic type.} Malabar: Krieger Pub. Company  


\bibitem{Hsu1} 
Hsu, E.~P. (1986). ``On excursions of reflecting Brownian motion.'' \emph{Trans. A.M.S.} 296 (1): 239--264.

\bibitem{Hsu2}
Hsu, E.~P. (2002).  ``Multiplicative functional for the heat equation on manifolds with boundary.'' \emph{Michigan Math. J.} 50 (2): 351--367.

\bibitem{IkeWat} 
Ikeda, N., Watanabe, S., (1981). \emph{Stochastic Differential Equations and Diffusion Processes.} Amsterdam: North Holland  

\bibitem{Jak1}
Jakubowski, A. (1997). ``A Non-Skorohod topology on the Skorohod space.'' \emph{Elec. J. Prob.} 2: 1--21.

\bibitem{Jak2}
Jakubowski, A. (2018). ``New characterizations of the S topology on the Skorokhod space.'' \emph{Elec. Comm. Prob.} 23: 1--16.

\bibitem{Kun} 
Kunita, H. (1997). \emph{Stochastic flows and stochastic differential equations.} Cambridge: Camb. Univ. Press  

\bibitem{LipRam} 
Lipschutz, D., Ramanan, K. (2019). ``Pathwise 
differentiability of reflected diffusions in convex polyhedral domains.'' \emph{Ann. I.H.P.} 55 (3 ): 1439--1476.

\bibitem{ManMas} 
Mandelbaum, A., Massey, W.~A. (1995). ``Strong approximations for time-dependent queues.'' \emph{Math. Oper. Res.} 20 (1): 33--64.

\bibitem{MitTay} 
Mitrea, D., Mitrea M., Mitrea I., Taylor, M. (1997). \emph{The Hodge-Laplacian: boundary value problems on Riemannian manifolds.} De Gruyter  

\bibitem{MeyZhe} 
Meyer, P.~A., Zheng, W.~A. . (1985). ``Tightness results for laws of diffusion processes application to stochastic mechanics.'' \emph{Ann. I.H.P.} 50 (2): 351--367.

\bibitem{Pil} 
Pilipenko, A.~Yu. (2007). ``Generalized differentiability with with respect to the initial data of a flow generated by a stochastic equation with reflection. '' \emph{Theory Prob. Math. Stat.} 75: 147--160.

\bibitem{RenWu} 
Ren, J., Wu, J. (2019). ``Probabilistic approach for nonlinear partial differential equations and stochastic partial differential equations with Neumann boundary conditions. '' \emph{J. Math. Anal. Appl.} 477: 1--40. 

\bibitem{RevYor}
Revuz, D., Yor, M. (1999). \emph{Continuous martingales and Brownian motion.} Springer 
\bibitem{Sco} 
Scott, C. (1995). ``$L^{p}$ theory of differential forms on manifolds.'' \emph{Trans. A.M.S.} 347 (6): 20
75--2096. 

\bibitem{Spi} 
Spivak, M. (1975). \emph{A comprehensive introduction to differential geometry. Vol. 5} Boston: Publish or Perish 
\bibitem{Tai} 
Taira, K. (2020) \emph{Boundary value problems and Markov processes.} Springer

\end{thebibliography}
\end{document}